\documentclass[11pt]{article}%
\usepackage{amsmath}
\usepackage{amsfonts}
\usepackage{amssymb}
\usepackage{graphicx}
\usepackage{a4wide}
\newtheorem{theorem}{Theorem}
\newtheorem{corollary}[theorem]{Corollary}

\newtheorem{lemma}[theorem]{Lemma}
\newtheorem{proposition}[theorem]{Proposition}
\newtheorem{remark}[theorem]{Remark}
\begin{document}

\title{New proofs for some results on spherically convex sets}
\author{Constantin Z\u{a}linescu\thanks{Octav Mayer Institute of Mathematics, Ia\c{s}i
Branch of Romanian Academy, Ia\c{s}i, Romania, and University
\textquotedblleft Alexandru Ioan Cuza" Ia\c{s}i, Romania; email:
\texttt{zalinesc@uaic.ro}.}}
\date{}
\maketitle

\begin{abstract}
In Guo and Peng's article [Spherically convex sets and spherically convex
functions, J. Convex Anal. 28 (2021), 103--122], one defines the notions of
spherical convex sets and functions on ``general curved surfaces'' in
$\mathbb{R}^{n}$ $(n\ge2)$, one studies several properties of these classes of
sets and functions, and one establishes analogues of Radon, Helly,
Carath\'{e}odory and Minkowski theorems for spherical convex sets, as well as
some properties of spherical convex functions which are analogous to those of
usual convex functions. In obtaining such results, the authors use an analytic
approach based on their definitions. Our aim in this note is to provide
simpler proofs for the results on spherical convex sets; our proofs are based
on some characterizations/representations of spherical convex sets by usual
convex sets in $\mathbb{R}^{n}$.

\end{abstract}

\section{Preliminaries}

As in \cite{GuoPen:21}, we consider $\mathbb{R}^{n}$ with $n\geq2$ endowed
with its usual inner product $\left\langle \cdot,\cdot\right\rangle $ and the
corresponding Euclidean norm $\left\Vert \cdot\right\Vert $, as well as a
continuous function $\Phi:\mathbb{R}^{n}\rightarrow\mathbb{R}_{+}
:=[0,\infty\lbrack$, such that $\Phi(tx)=t\Phi(x)$ for all $x\in\mathbb{R}%
^{n}$, $t\in\mathbb{R}_{+}$, and $\Phi(x)=0$ if and only if $x=o$, where $o$
denotes the zero element of $\mathbb{R}^{n}$. The set $\mathbb{S}
:=\mathbb{S}_{\Phi}:= \{x\in\mathbb{R}^{n} \mid\Phi(x)=1\}$ is called $\Phi
$-sphere, or simply sphere when there is no danger of confusion. We also
consider the mapping $\rho:=\rho_{\Phi}:\mathbb{R}^{n}\rightarrow
\{o\}\cup\mathbb{S}$ defined by $\rho(x):=[\Phi(x)]^{-1}x$ if $x\neq o$ and
$\rho(o):=o$; clearly, $\rho(tx)=\rho(x)$ for all $x\in\mathbb{R}^{n}$ and
$t\in\mathbb{P}:={}]0,\infty\lbrack$. On $\mathbb{S}$ we consider the trace of
the (usual) topology of $\mathbb{R}^{n}$. It follows that $\mathbb{S}$ is
compact. Clearly, when $\Phi:=\left\Vert \cdot\right\Vert $, $\mathbb{S}%
_{\Phi}=\mathbb{S}^{n-1}$ is the $n-1$ dimensional (unit) sphere.

Having $x_{1},...,x_{k}\in\mathbb{S}$ and $\lambda_{1},...,\lambda_{k}%
\in\mathbb{R}_{+}$ $(k\geq1)$ such that $\sum_{i=1}^{k}\lambda_{i}=1$,
\[
\text{(s)}%
%TCIMACRO{\tsum \nolimits_{i=1}^{k}}%
%BeginExpansion
{\textstyle\sum\nolimits_{i=1}^{k}}
%EndExpansion
\lambda_{i}x_{i}:=\text{(s}_{\Phi}\text{)}%
%TCIMACRO{\tsum \nolimits_{i=1}^{k}}%
%BeginExpansion
{\textstyle\sum\nolimits_{i=1}^{k}}
%EndExpansion
\lambda_{i}x_{i}:=\rho\big(%
%TCIMACRO{\tsum \nolimits_{i=1}^{k}}%
%BeginExpansion
{\textstyle\sum\nolimits_{i=1}^{k}}
%EndExpansion
\lambda_{i}x_{i}\big)
\]
is called a \emph{$\Phi$-spherical convex} (\emph{s-convex} for short)
\emph{combination} of $x_{1},...,x_{k};$ when $x,y\in\mathbb{S}$ and
$\lambda\in\lbrack0,1]$, $\rho(\lambda x+(1-\lambda)y)$ is written $\lambda
x+_{s}(1-\lambda)y$. The nonempty subset $S\subseteq\mathbb{S}$ is called
\emph{s-convex} if $\lambda x+_{s}(1-\lambda)y\in S$ for all $x,y\in S$ and
$\lambda\in\lbrack0,1];$ of course, we may take $\lambda\in{}]0,1[$ in the
preceding definition.

\smallskip We denote by $\mathcal{C}_{\text{s}}$ the class of the nonempty
subsets of $\mathbb{S}$ which are s-convex.

\smallskip For $\emptyset\neq A\subseteq\mathbb{R}^{n}$ we denote by
$\operatorname*{span}A$, $\operatorname*{aff}A$, $\operatorname*{conv}A$,
$\overline{\operatorname*{conv}}A$, $\operatorname*{cone}A$, $\overline
{\operatorname*{cone}}A$, $\operatorname*{cl}A$, $\operatorname*{bd}A$ the
linear hull, affine hull, the convex hull, the closed convex hull, the convex
conic hull, the closed convex conic hull, the closure and the boundary of $A$,
respectively; recall that a cone $K\subseteq\mathbb{R}^{n}$ is pointed if
$K\cap(-K)=\{o\}$. Moreover, $\Gamma A:=\{\gamma a\mid\gamma\in\Gamma$, $a\in
A\}$ and $A+B:=\{a+b\mid a\in A$, $b\in B\}$ for $\emptyset\neq\Gamma
\subseteq\mathbb{R}$ and $\emptyset\neq B\subseteq\mathbb{R}^{n}$. Clearly,
$\operatorname*{cone}A=\mathbb{R}_{+}\cdot\operatorname*{conv}%
A=\operatorname*{conv}(\mathbb{R}_{+}A)$.

\begin{proposition}
\label{p1}Let $\emptyset\neq S\subseteq\mathbb{S}$. The following assertions
are equivalent:

\emph{(i)} $S$ is s-convex;

\emph{(ii)} $\forall x_{1},x_{2}\in\mathbb{P}S:x_{1}+x_{2}\in\mathbb{P}S$;

\emph{(iii)} $\mathbb{P}S$ is convex;

\emph{(iv)} $\mathbb{R}_{+}S$ is a pointed, convex cone;

\emph{(v)} there exists a pointed convex cone $K\subseteq\mathbb{R}^{n}$ such
that $S=K\cap\mathbb{S}$;

\emph{(vi)} there exists a convex set $C\subseteq\mathbb{R}^{n}$ such that
$S=\rho(C)$;

\emph{(vii)} there exists a bounded convex set $C\subseteq\mathbb{R}^{n}$ such
that $S=\rho(C)$;

\emph{(viii)} $S=\rho\left(  \operatorname*{conv}S\right)  $.

Assume, moreover, that $\Phi$ is convex. Then \emph{(i)}~$\Leftrightarrow$
\emph{(ix)}, where

\emph{(ix)} $]0,1]S$ is convex.
\end{proposition}

Proof. The implications (ix) $\Rightarrow$ (vii) $\Rightarrow$ (vi), (viii)
$\Rightarrow$ (vii) and (iv) $\Rightarrow$ (v) are obvious.

(i) $\Rightarrow$ (ii) Take $x_{1},x_{2}\in\mathbb{P}S$; then there exist
$u_{i}\in S$ and $t_{i}\in\mathbb{P}$ such that $x_{i}=t_{i}u_{i}$ for
$i=1,2$. Then $x_{1}+x_{2}=(t_{1}+t_{2})u$, where $u:=\frac{t_{1}}{t_{1}%
+t_{2}}u_{1}+\frac{t_{2}}{t_{1}+t_{2}}u_{2}$. By the s-convexity of $S$,
$\rho(u)\in S$, and so $u=\Phi(u)\cdot\rho(u)\in\mathbb{P}S$.

(ii) $\Rightarrow$ (iii) It is obvious that $tx\in\mathbb{P}S$ if
$t\in\mathbb{P}$ and $x\in\mathbb{P}S$. Take $x_{1},x_{2}\in\mathbb{P}S$ and
$\lambda\in{}]0,1[$. Then $\lambda x_{1}$, $(1-\lambda)x_{2}\in\mathbb{P}S$
and so $\lambda x_{1}+(1-\lambda)x_{2}\in\mathbb{P}S$. Hence $\mathbb{P}S$ is convex.

(iii) $\Rightarrow$ (iv) The convexity of $\mathbb{R}_{+}S$ is (almost)
obvious. Let $\pm x\in(\mathbb{R}_{+}S)\setminus\{o\}$; then $\pm
x\in\mathbb{P}S$, whence the contradiction $o=\tfrac{1}{2}x+\tfrac{1}%
{2}(-x)\in\mathbb{P}S$. Therefore, $\mathbb{R}_{+}S$ is a pointed convex cone.

(v)\ $\Rightarrow$ (i) Take $x_{1},x_{2}\in S$ $(\subseteq K\setminus\{o\})$
and $\lambda\in{}]0,1[$. Then $u:=\lambda x_{1}+(1-\lambda)x_{2}\in K$. We
claim that $u\neq o$; in the contrary case $K\setminus\{o\}\ni x:=\lambda
x_{1}=(\lambda-1)x_{2}=-x\in K$, whence the contradiction $x=o$. It follows
that $\Phi(u)>0$, and so $\rho(u)=\frac{u}{\Phi(u)}\in K\cap\mathbb{S}=S$.
Hence $S$ is s-convex.

(iii) $\Rightarrow$ (viii) Because $\mathbb{P}S$ is convex and $S\subseteq
\mathbb{P}S$, we have $(S\subseteq)$ $\operatorname*{conv}S\subseteq
\mathbb{P}S$, and so $S=\rho(S)\subseteq\rho(\operatorname*{conv}%
S)\subseteq\rho(\mathbb{P}S)=\rho(S)$, whence $S=\rho(\operatorname*{conv}S)$.

(vi) $\Rightarrow$ (iii) We have that $\mathbb{P}S=\mathbb{P}\cdot
\rho(C)=\mathbb{P}C$, and so $\mathbb{P}S$ is convex because $C$ is so.

(i) $\Rightarrow$ (ix) Assume that $\Phi$ is convex. Take $x_{1},x_{2}\in
{}]0,1]S$; then there exist $u_{i}\in S$ and $t_{i}\in{}]0,1]$ such that
$x_{i}=t_{i}u_{i}$ for $i=1,2$. Consider $\lambda\in{}]0,1[$. Then $x:=\lambda
x_{1}+(1-\lambda)x_{2}=\gamma y$, where $\gamma:=\lambda t_{1}+(1-\lambda
)t_{2}\in{}]0,1]$ and $y:=\mu u_{1}+(1-\mu)u_{2}$ with $\mu:=\gamma
^{-1}\lambda t_{1}\in{}]0,1[$. Because $S$ is s-convex, we have that
$\rho(y)\in S$, while because $\Phi$ is convex one has $(0<)$ $\Phi(y)\leq1$,
and so $\gamma\Phi(y)\in{}]0,1]$. It follows that $x=\gamma\Phi(y)\rho(y)\in
{}]0,1]S$. \hfill$\square$

\medskip The equivalence of assertions (i)--(iv) from Proposition \ref{p1} is
established in \cite[Prop.\ 2.7]{GuoPen:21}. However, the equivalence
(i)\ $\Leftrightarrow$ (iv) in the case $\Phi:=\left\Vert \cdot\right\Vert $
seems to be known for a long time (see e.g.\ \cite[p.\ 158]{DaGrKl:63} where
one says: \textquotedblleft Due to the close and obvious relationship between
spherical convexity in $S^{n}$ and convex cones in $R^{n+1}$, many results can
be interpreted in both the spherical and the Euclidean setting"); this
equivalence is mentioned without proof in \cite[p.\ 11]{Bak:94} and with proof
in \cite[Prop.\ 2]{FeIuNe:13}.

\begin{remark}\emph{
\label{rem1}Clearly, the convex sets $C$ from (vi) and (vii) do not
contain the origin $o$ of $\mathbb{R}^{n}$, and so
$o\notin\operatorname*{conv}S$ if $S\in\mathcal{C}_{\text{s}}$.
Moreover, the equivalences (i)~$\Leftrightarrow$
(v)~$\Leftrightarrow$ (vi)~$\Leftrightarrow$ (vii) describe
completely which are the s-convex subsets of $\mathbb{S}$; more
precisely,}
\begin{align}
\mathcal{C}_{\text{s}} &  =\big\{K\cap\mathbb{S}\mid\mathbb{R}^{n}\supseteq
K=\operatorname*{cone}K\neq\{o\},\ K\cap(-K)=\{o\}\big\}\nonumber\\
&  =\{\rho(C)\mid\mathbb{R}^{n}\setminus\{o\}\supseteq C=\operatorname*{conv}%
C\neq\emptyset\}\nonumber\\
&  =\{\rho(C)\mid\mathbb{R}^{n}\setminus\{o\}\supseteq C=\operatorname*{conv}%
C\neq\emptyset,\ C\ \text{\rm is bounded}\}.\label{r-Cs}%
\end{align}

\end{remark}

From Proposition \ref{p1} one obtains rapidly several characterizations of the
closed s-convex subsets of $\mathbb{S}$; just observe that $\mathbb{R}_{+}B$
is closed whenever $B$ is a bounded and closed subset of $\mathbb{R}^{n}$ such
that $o\notin B$, the closed subsets of $\mathbb{S}$ are compact and $\rho(E)$
is compact for every compact subset $E$ of $\mathbb{R}^{n}\setminus\{o\}$
(because $\rho|_{\mathbb{R}^{n}\setminus\{o\}}$ is continuous).

\begin{corollary}
\label{c-p1}Let $\emptyset\neq S\subseteq\mathbb{S}$. The following assertions
are equivalent:

\emph{(i)} $S$ is closed and s-convex;

\emph{(ii)} $\mathbb{R}_{+}S$ is a pointed, closed and convex cone;

\emph{(iii)} there exists a pointed, closed and convex cone $K\subseteq
\mathbb{R}^{n}$ such that $S=K\cap\mathbb{S}$;

\emph{(iv)} $S=\rho\left(  \overline{\operatorname*{conv}}S\right)  $;

\emph{(v)} there exists a compact and convex set $C\subseteq\mathbb{R}^{n}$
such that $S=\rho\left(  C\right)  $.
\end{corollary}

We denote by $\mathcal{C}_{\text{sc}}$ the class of the closed and s-convex
subsets of $\mathbb{S}$; hence%
\begin{align}
\mathcal{C}_{\text{sc}}  &  =\big\{K\cap\mathbb{S}\mid\mathbb{R}^{n}\supseteq
K=\overline{\operatorname*{cone}}K\neq\{o\},\ K\cap(-K)=\{o\}\big\}\nonumber\\
&  =\{\rho(C)\mid\mathbb{R}^{n}\setminus\{o\}\supseteq C=\overline
{\operatorname*{conv}}C\neq\emptyset,\ C\ \text{is bounded}\}. \label{r-Csc}%
\end{align}

\section{Proofs}

In the sequel we provide proofs of the main results on spherical convex sets
established in \cite{GuoPen:21}; they are based on Proposition \ref{p1}, as
well as on other results stated in the sequel.

\begin{theorem}
\label{cor-GPt3.1}\emph{(\cite[Th.\ 3.1]{GuoPen:21})} Let $C\subseteq
\mathbb{S}$ be an s-convex set. Then for $(x_{i})_{i=1}^{k}\subseteq C$ and
$(\lambda_{i})_{i=1}^{k}\subseteq\lbrack0,1]$ with $\sum_{i=1}^{k}\lambda
_{i}=1$, one has \emph{(s)}$\sum_{i=1}^{k}\lambda_{i}x_{i}\in C$. In
particular, $o\notin\operatorname*{conv}C$.
\end{theorem}

Proof. Take $(x_{i})_{i=1}^{k}\subseteq C$ and $(\lambda_{i})_{i=1}%
^{k}\subseteq\mathbb{R}_{+}$ with $\sum_{i=1}^{k}\lambda_{i}=1$. Then
$x:=\sum_{i=1}^{k}\lambda_{i}x_{i}\in\operatorname*{conv}C$. Using the
implication (i)~$\Rightarrow$ (viii) of Proposition \ref{p1} we obtain that
(s)$\sum_{i=1}^{k}\lambda_{i}x_{i}=\rho(x)\in\rho(\operatorname*{conv}C)=C$.
From Remark \ref{rem1} one has $o\notin\operatorname*{conv}C$. \hfill
\hfill$\square$

\begin{theorem}
\emph{(\cite[Th.\ 3.2]{GuoPen:21})} \label{cor-GPt3.2} Let $C\subseteq
\mathbb{S}$ be an s-convex set. Then:

\emph{(i)} There is $u_{0}\in\mathbb{S}^{n-1}$ such that $C\subseteq
\overline{H}_{u_{0}}^{+}\cap\mathbb{S}$, where $\overline{H}_{u_{0}}%
^{+}:=[\left\langle \cdot,u_{0}\right\rangle \geq0]:=\{x\in\mathbb{R}^{n}%
\mid\left\langle x,u_{0}\right\rangle \geq0\}$, and so $C$ is contained in a
closed $\Phi$-hemisphere.

\emph{(ii)} If further $C$ is closed, there are $u_{0}\in\mathbb{S}^{n-1}$ and
$\alpha>0$ such that $C\subseteq\overline{H}_{u_{0},\alpha}^{+}\cap\mathbb{S}%
$, where $\overline{H}_{u_{0},\alpha}^{+}:=[\left\langle \cdot,u_{0}%
\right\rangle \geq\alpha]$ $(\subseteq\lbrack\left\langle \cdot,u_{0}%
\right\rangle >0])$. In particular, $C$ is contained in an open $\Phi$-hemisphere.
\end{theorem}

Proof. (i) By Remark \ref{rem1} we have that $o\notin\operatorname*{conv}C$,
and so there exists $u_{0}\in\mathbb{S}^{n-1}$ such that $\left\langle
x,u_{0}\right\rangle \geq0$ for every $x\in\operatorname*{conv}C$; hence
$C\subseteq\operatorname*{conv}C\subseteq\overline{H}_{u_{0}}^{+}$, and so
$C\subseteq\overline{H}_{u_{0}}^{+}\cap\mathbb{S}$.

(ii) Because $C$ is a closed subset of the compact set $\mathbb{S}$,
$\operatorname*{conv}C$ is compact. Because $o\notin\operatorname*{conv}C$,
there exist $u_{0}\in\mathbb{S}^{n-1}$ and $\alpha\in\mathbb{R}$ such that
$\left\langle x,u_{0}\right\rangle \geq\alpha>0$ for every $x\in
\operatorname*{conv}C$. As above we get $C\subseteq\overline{H}_{u_{0},\alpha
}^{+}\cap\mathbb{S}$. \hfill$\square$

\medskip The next result is surely known. In the case $\dim X<\infty$,
assertions (i) and (ii) follow by Corollaries 17.1.1 and 17.1.2 from
\cite{Roc:70}, respectively.

\begin{lemma}
\label{l-pg1}Consider $X$ a real linear space and $\emptyset\neq A\subseteq
X$. Then:

\emph{(i)} $\operatorname*{conv}A=\big\{\sum_{k=1}^{m}\alpha_{k}x_{k}\mid
m\geq1$, $(\alpha_{k})_{_{k=1}}^{m}\subseteq\mathbb{P}$, $(x_{k})_{_{k=1}}%
^{m}\subseteq A$ a.i., $\sum_{k=1}^{m}\alpha_{k}=1\big\}$;

\emph{(ii)} $\mathbb{P}\cdot\operatorname*{conv}A=\big\{\sum_{k=1}^{m}%
\alpha_{k}x_{k}\mid m\geq1$, $(\alpha_{k})_{_{k=1}}^{m}\subseteq\mathbb{P}$,
$(x_{k})_{_{k=1}}^{m}\subseteq A$ l.i.$\big\}$.
\end{lemma}

In the text above (and later on), ``a.i." and ``l.i." are abbreviations for
\textquotedblleft affinely independent\textquotedblright\ and
\textquotedblleft linearly independent\textquotedblright, respectively.

\begin{proposition}
\label{p3}Let $\emptyset\neq C\subseteq\mathbb{R}^{n}\setminus\{o\}$. Consider
the following assertions:

\emph{(i)} there exists $u_{0}\in\mathbb{R}^{n}$ such that
$\left\langle x,u_{0}\right\rangle >0$ for all $x\in C$;

\emph{(ii)} for every a.i.\ family $(x_{i})_{_{i=1}}^{k}\subseteq C$ and
$(\alpha_{i})_{_{i=1}}^{k}\subseteq\mathbb{P}$ with $k\geq1$ one has
$\sum\nolimits_{i=1}^{k}\alpha_{i}x_{i}\neq o$;

\emph{(iii)} for every a.i.\ family $(x_{i})_{_{i=1}}^{k}\subseteq C$ with
$k\geq1$ one has $\mathbb{P}\cdot\operatorname*{conv}\{x_{i}\mid i\in
\overline{1,k}\}\neq\operatorname*{span}\{x_{i}\mid i\in\overline{1,k}\}$.

Then \emph{(i)}~$\Rightarrow$ \emph{(ii)}~$\Leftrightarrow$ \emph{(iii)};
moreover, if $C$ is compact, then \emph{(ii)}~$\Rightarrow$ \emph{(i)}.
\end{proposition}

Proof. (i)~$\Rightarrow$ (ii) Take an a.i.\ family $(x_{i})_{_{i=1}}%
^{k}\subseteq C$ and $(\alpha_{i})_{_{i=1}}^{k}\subseteq\mathbb{P}$ with
$k\geq1$. Then $\big\langle\sum\nolimits_{i=1}^{k}\alpha_{i}x_{i}%
,u_{0}\big\rangle=\sum\nolimits_{i=1}^{k}\alpha_{i}\left\langle x_{i}%
,u_{0}\right\rangle >0$, and so $\sum\nolimits_{i=1}^{k}\alpha_{i}x_{i}\neq o$.

(ii)~$\Rightarrow$ (i) Assume that $C$ is compact. Using Lemma
\ref{l-pg1} (i), one obtains that $o\notin\operatorname*{conv}C$.
Because $C$ is compact, $\operatorname*{conv}C$ is closed; using a
separation theorem, there exists $u_{0}\in\mathbb{R}^{n}$ such that
$\langle x,u_{0}\rangle>0$ for all $x\in\operatorname*{conv}C$, that
is (i) holds.

(ii)~$\Leftrightarrow$ (iii) In fact, for a nontrivial real linear space $X$
and a nonempty convex set $A\subseteq X$ the following equivalences hold:%
\begin{align*}
0\in\operatorname*{icr}A  &  \Longleftrightarrow\mathbb{R}_{+}A\text{ is a
linear space}\Longleftrightarrow\mathbb{R}_{+}A=\mathbb{R}_{+}A-\mathbb{R}%
_{+}A\Longleftrightarrow\mathbb{R}_{+}A=\operatorname*{span}A\\
&  \Longleftrightarrow\mathbb{P}A\text{ is a linear space}\Longleftrightarrow
\mathbb{P}A=\mathbb{P}A-\mathbb{P}A\Longleftrightarrow\mathbb{P}%
A=\operatorname*{span}A,
\end{align*}
where $\operatorname*{icr}E:=\{x\in X\mid\forall x^{\prime}\in
\operatorname*{aff}E$, $\exists\delta>0$, $\forall t\in\lbrack0,\delta
]:(1-t)x+tx^{\prime}\in E\}$ is the relative algebraic interior (or intrinsic
core) of $\emptyset\neq E\subseteq X$. Observing that for $\emptyset\neq
E\subseteq X$ one has $\operatorname*{span}E=\operatorname*{span}%
(\operatorname*{conv}E)$, from the above equivalences applied for
$A:=\operatorname*{conv}\{x_{i}\mid i\in\overline{1,k}\}$, one has that (iii)
is equivalent to $0\notin\operatorname*{icr}(\operatorname*{conv}\{x_{i}\mid
i\in\overline{1,k}\})$ for every a.i.\ family $(x_{i})_{_{i=1}}^{k}\subseteq
C$. Hence (ii)~$\Leftrightarrow$ (iii) because for $(x_{i})_{i\in
\overline{1,k}}\subseteq X$ an a.i.\ family one has $\operatorname*{icr}%
(\operatorname*{conv}\{x_{i}\mid i\in\overline{1,k}\}=\big\{\sum
\nolimits_{i=1}^{k}\alpha_{i}x_{i}\mid(\alpha_{i})_{_{i=1}}^{k}\subseteq
\mathbb{P}$, $\sum\nolimits_{i=1}^{k}\alpha_{i}=1\big\}$. \hfill$\square$

\medskip Observe that \cite[Th.\ 3.3]{GuoPen:21} asserts that $\rceil
$(i)~$\Leftrightarrow$ $\rceil$(ii)~$\Leftrightarrow$ $\rceil$(iii) for
$C\subseteq\mathbb{S}^{n-1}$ a compact set, where (i)--(iii) are the
assertions from Proposition \ref{p3}.

As in \cite{GuoPen:21}, one says that $\emptyset\neq S\subseteq\mathbb{S}$ is
\emph{hull-addible} if $o\notin\operatorname*{conv}S;$ moreover, for
$S\subseteq\mathbb{S}$ hull-addible one sets $\operatorname*{sco}%
S:=\rho(\operatorname*{conv}S)$.

\begin{theorem}
\emph{(\cite[Th.\ 3.8]{GuoPen:21})} \label{cor-GPt3.6} Let $S\subseteq
\mathbb{S}$ be a hull-addible set. Then

\emph{(i)} $\operatorname*{sco}S$ is s-convex.

\emph{(ii)} $\operatorname*{sco}S=%
%TCIMACRO{\tbigcap }%
%BeginExpansion
{\textstyle\bigcap}
%EndExpansion
\{C\in\mathcal{C}_{\text{s}}\mid S\subseteq C\}$, i.e. $\operatorname*{sco}S$
is the smallest s-convex set containing $S$.

\emph{(iii)} $\operatorname*{sco}S=\mathbb{S}\cap\operatorname*{cone}S$.
\end{theorem}

Proof. (i) Taking into account (\ref{r-Cs}), one has $\operatorname*{sco}%
S\in\mathcal{C}_{\text{s}}$ because $\operatorname*{sco}S:=\rho
(\operatorname*{conv}S)$.

(ii) The inclusion $\supseteq$ follows from (i); the inclusion $\subseteq$ is
true because for $C\in\mathcal{C}_{\text{s}}$ with $S\subseteq C$
$(\subseteq\mathbb{S})$ we have that $\operatorname*{conv}S\subseteq
\operatorname*{conv}C$, and so $\operatorname*{sco}S:=\rho
(\operatorname*{conv}S)\subseteq\rho(\operatorname*{conv}C)=C$, the last
equality being given by the implication (i)~$\Rightarrow$ (viii) of
Proposition \ref{p1}.

(iii) Having in view that
\[
\mathbb{P}\cdot\operatorname*{sco}S=\mathbb{P}\cdot\mathbb{\rho(}%
\operatorname*{conv}S)=\mathbb{P}\cdot\operatorname*{conv}S=\left(
\mathbb{R}_{+}\cdot\operatorname*{conv}S\right)  \setminus
\{o\}=(\operatorname*{cone}S)\setminus\{o\},
\]
we obtain that $\mathbb{S}\cap\operatorname*{cone}S=\mathbb{S}\cap
(\mathbb{P}\cdot\operatorname*{sco}S)=\operatorname*{sco}S$. \hfill
\hfill$\square$

\begin{proposition}
\emph{(\cite[Rem.\ 3.11]{GuoPen:21})} \label{cor-GPp3.8}Let $S\subseteq
\mathbb{S}$ be a nonempty set. Then there is a closed s-convex set containing
$S$ if and only if there exist $u\in\mathbb{S}^{n-1}$ and $\alpha>0$ such that
$S\subseteq\mathbb{S}\cap\overline{H}_{u,\alpha}^{+}$.
\end{proposition}

Proof. Assume that $S\subseteq\overline{H}_{u,\alpha}^{+}$ for some
$u\in\mathbb{S}^{n-1}$ and $\alpha>0$. Then $C:=\overline{\operatorname*{conv}%
}S\subseteq\overline{H}_{u,\alpha}^{+}$ is a compact convex set with $o\notin
C\supseteq S$. Using the implication (v)~$\Rightarrow$ (i) of Corollary
\ref{c-p1} we have that $S^{\prime}:=\rho(C)\in\mathcal{C}_{\text{sc}}$ and
$S\subseteq S^{\prime}$. Assume now that $S\subseteq S^{\prime}$ with
$S^{\prime}\in\mathcal{C}_{\text{sc}}$. By Theorem \ref{cor-GPt3.2} (ii),
$(S\subseteq)$ $S^{\prime}\subseteq\mathbb{S}\cap\overline{H}_{u,\alpha}^{+}$
for some $u\in\mathbb{S}^{n-1}$ and $\alpha>0$. \hfill$\square$

\medskip Notice that Proposition \ref{cor-GPp3.8} is established in
\cite[Prop.\ 3.10]{GuoPen:21} for $\Phi$ convex, and stated without proof in
\cite[Rem.\ 3.11]{GuoPen:21}.

\begin{theorem}
\emph{(\cite[Th.\ 3.12]{GuoPen:21})} \label{cor-GPt3.10}Let $\emptyset\neq
S\subseteq\mathbb{S}\cap\overline{H}_{u,\alpha}^{+}$ for some $u\in
\mathbb{S}^{n-1}$ and $\alpha>0$. If $\Phi$ is convex, then%
\begin{equation}
\overline{\operatorname*{sco}}S:=\operatorname*{cl}(\operatorname*{sco}S)=%
%TCIMACRO{\tbigcap }%
%BeginExpansion
{\textstyle\bigcap}
%EndExpansion
\{C\in\mathcal{C}_{\text{sc}}\mid S\subseteq C\}. \label{r-scob}%
\end{equation}

Consequently, $\overline{\operatorname*{sco}}S$ is the smallest closed and
s-convex subset of $\mathbb{S}$ containing $S$.
\end{theorem}

Proof. By Corollary \ref{cor-GPp3.8}, there exists $C_{0}\in\mathcal{C}%
_{\text{sc}}$ such that $S\subseteq C_{0}$; it follows that
$\operatorname*{cl}(\operatorname*{sco}S)\subseteq C_{0}$. Then $(\emptyset
\neq)$ $S\subseteq\Sigma:=%
%TCIMACRO{\tbigcap }%
%BeginExpansion
{\textstyle\bigcap}
%EndExpansion
\{C\in\mathcal{C}_{\text{sc}}\mid S\subseteq C\}\in\mathcal{C}_{\text{sc}}$,
and so $\operatorname*{cl}(\operatorname*{sco}S)\subseteq\Sigma$. On the other
hand, it is known (and easy to prove) that $\operatorname*{cl}(\mathbb{R}%
_{+}B)=\mathbb{R}_{+}\cdot\operatorname*{cl}B$ whenever $B\subseteq
\mathbb{R}^{n}$ is bounded and $0\notin\operatorname*{cl}B$. (In fact, this
assertion is true for $\mathbb{R}^{n}$ replaced by an arbitrary Hausdorff
topological vector space.) Since $\operatorname*{cl}(\operatorname*{sco}%
S)\subseteq C_{0}$, one has
\[
K:=\mathbb{R}_{+}\cdot\operatorname*{cl}(\operatorname*{sco}%
S)=\operatorname*{cl}(\mathbb{R}_{+}\cdot\operatorname*{sco}S)\subseteq
\mathbb{R}_{+}C_{0}.
\]
Consequently, $K$ is pointed (since $\mathbb{R}_{+}C_{0}$ is so), convex
(since $\mathbb{R}_{+}\cdot\operatorname*{sco}S$ is so) and closed, whence
$(S\subseteq)$ $\operatorname*{cl}(\operatorname*{sco}S)\in\mathcal{C}%
_{\text{sc}}$. Therefore, $\Sigma\subseteq\operatorname*{cl}%
(\operatorname*{sco}S)$, and so (\ref{r-scob}) holds. \hfill$\square$

\medskip Notice that Theorem \ref{cor-GPt3.10} is established in
\cite[Th.\ 3.12]{GuoPen:21} for $\Phi$ convex.

\smallskip The following variant of Radon's theorem holds in $\mathbb{R}^{n}$.

\begin{lemma}
\label{l-GPt4.1} Let $A:=\{u_{i}\mid i\in I\}\subseteq\mathbb{R}^{n}$ be such
that $\operatorname*{card}I\geq n+1$ and $o\notin\operatorname*{conv}A$. Then
there exist two nonempty sets $I_{1},I_{2}\subseteq I$ such that $I_{1}\cap
I_{2}=\emptyset$, $I_{1}\cup I_{2}=I$ and $\operatorname*{conv}A_{1}%
\cap\left(  ]0,1]\operatorname*{conv}A_{2}\right)  \neq\emptyset$, where
$A_{k}:=\{u_{i}\mid i\in I_{k}\}$ for $k\in\{1,2\}$.
\end{lemma}

Proof. Take $i^{\prime}\notin I$, and set $u_{i^{\prime}}:=o$, $I^{\prime
}:=I\cup\{i^{\prime}\}$ and $A^{\prime}:=A\cup\{u_{i^{\prime}}\}$; then
$\operatorname*{card}I^{\prime}\geq n+2$. By Radon's theorem (see
\cite[Th.\ 1.1.5]{Sch:93}), there exist $\emptyset\neq I_{1}^{\prime}%
,I_{2}^{\prime}\subseteq I^{\prime}$ such that $I_{1}^{\prime}\cap
I_{2}^{\prime}=\emptyset$, $I_{1}^{\prime}\cup I_{2}^{\prime}=I^{\prime}$ and
$\operatorname*{conv}A_{1}^{\prime}\cap\operatorname*{conv}A_{2}^{\prime}%
\neq\emptyset$, where $A_{k}^{\prime}:=\{u_{i}\mid i\in I_{k}^{\prime}\}$ for
$k\in\{1,2\}$. Clearly $I_{k}^{\prime}\neq\{i^{\prime}\}$ for $k\in\{1,2\}$;
otherwise we would have $o\in\operatorname*{conv}A$. Let $i^{\prime}\in
I_{2}^{\prime}$ and set $I_{1}:=I_{1}^{\prime}$, $I_{2}:=I_{2}^{\prime
}\setminus\{i^{\prime}\}$. Then $\operatorname*{conv}A_{1}^{\prime
}=\operatorname*{conv}A_{1}$ and $\operatorname*{conv}A_{2}^{\prime
}=\operatorname*{conv}\left(  A_{2}\cup\{o\}\right)
=[0,1]\operatorname*{conv}A_{2}$, whence the conclusion follows because
$o\notin\operatorname*{conv}A\supseteq\operatorname*{conv}A_{1}\cup
\operatorname*{conv}A_{2}$. \hfill$\square$

\begin{theorem}
\emph{(\cite[Th.\ 4.1]{GuoPen:21}, Radon-type Theorem)} Let $S\subseteq
\mathbb{S}$ be such that $\operatorname*{card}S\geq n+1$ and $o\notin
\operatorname*{conv}S$. Then there exist two nonempty sets $S_{1}%
,S_{2}\subseteq S$ such that $S_{1}\cap S_{2}=\emptyset$, $S_{1}\cup S_{2}=S$
and $\operatorname*{sco}S_{1}\cap\operatorname*{sco}S_{2}\neq\emptyset$.
\end{theorem}

Proof. Using the preceding lemma there exist two nonempty sets $S_{1}%
,S_{2}\subseteq S$ such that $S_{1}\cap S_{2}=\emptyset$, $S_{1}\cup S_{2}=S$
and $\operatorname*{conv}S_{1}\cap\left(  ]0,1]\operatorname*{conv}%
S_{2}\right)  \neq\emptyset$. Because $\operatorname*{sco}S_{1}=\rho
(\operatorname*{conv}S_{1})$ and $\operatorname*{sco}S_{2}=\rho
(\operatorname*{conv}S_{2})=\rho\left(  ]0,1]\operatorname*{conv}S_{2}\right)
$, the conclusion follows. \hfill$\square$

\medskip

The following result is surely known (see \cite[p.\ 713]{SarSas:78} for a
hint). We provide a simple proof for readers convenience.{}

\begin{lemma}
\label{l-ss}Let $(u_{i})_{i\in I}\subseteq\mathbb{R}^{n}$ be an a.i.\ family,
where $I:=\overline{1,n+1}$. Assume that $o\in\operatorname*{core}%
(\operatorname*{conv}\{u_{i}\mid i\in I\})$, or, equivalently, there exists
$(\overline{\lambda}_{l})_{l\in I}\subseteq\mathbb{P}$ such that
$\sum\nolimits_{l\in I}\overline{\lambda}_{l}=1$ and $\sum\nolimits_{l\in
I}\overline{\lambda}_{l}u_{l}=o$. For each $i\in I$ set
\begin{equation}
Q_{i}:=\left\{
%TCIMACRO{\tsum \nolimits_{l\in I}}%
%BeginExpansion
{\textstyle\sum\nolimits_{l\in I}}
%EndExpansion
\lambda_{l}u_{l}\mid(\lambda_{l})_{l\in I}\subseteq\mathbb{R}_{+},\text{
}\lambda_{i}=0\right\}  . \label{r-pg7}%
\end{equation}

Then: \emph{(a)}$~Q_{i}$ is a pointed convex cone for each $i\in I$;
\emph{(b)}~for every $\emptyset\neq J\subseteq I$ one has
\begin{equation}%
%TCIMACRO{\tbigcap \nolimits_{j\in J}}%
%BeginExpansion
{\textstyle\bigcap\nolimits_{j\in J}}
%EndExpansion
Q_{j}=\left\{
%TCIMACRO{\tsum \nolimits_{l\in I}}%
%BeginExpansion
{\textstyle\sum\nolimits_{l\in I}}
%EndExpansion
\lambda_{l}u_{l}\mid(\lambda_{l})_{l\in I}\subseteq\mathbb{R}_{+},\text{
}\lambda_{j}=0\ \forall j\in J\right\}  ; \label{r-pg4}%
\end{equation}
\emph{(c)}~$\bigcap\nolimits_{i\in I\setminus\{j\}}Q_{i}=\mathbb{R}_{+}u_{j}$
for every $j\in I$ and $\bigcap\nolimits_{i\in I}Q_{i}=\{0\}$; \emph{(d)}%
~$\bigcup\nolimits_{i\in I}Q_{i}=\mathbb{R}^{n}$.
\end{lemma}

Proof. Let us observe that having $x\in X:=\mathbb{R}^{n}$ such that
$x=\sum\nolimits_{l\in I}\lambda_{l}u_{l}=\sum\nolimits_{l\in I}\lambda
_{l}^{\prime}u_{l}$ with $\lambda_{l},\lambda_{l}^{\prime}\geq0$ for $l\in I$
and $\lambda_{i}=\lambda_{j}^{\prime}=0$ for some $i,j\in I$, we have that
$\lambda_{l}=\lambda_{l}^{\prime}$ for all $l\in I$. Indeed, setting
$\gamma:=\sum\nolimits_{l\in I}(\lambda_{l}-\lambda_{l}^{\prime})$, clearly
$\sum\nolimits_{l\in I}(\lambda_{l}-\lambda_{l}^{\prime}-\gamma\overline
{\lambda}_{l})=0$ and $\sum\nolimits_{l\in I}(\lambda_{l}-\lambda_{l}^{\prime
}-\gamma\overline{\lambda}_{l})u_{l}=o$, and so $\lambda_{l}-\lambda
_{l}^{\prime}=\gamma\overline{\lambda}_{l}$ for every $l\in I$. Taking
$l\in\{i,j\}$, we get $(0\geq)$ $-\lambda_{i}^{\prime}=\gamma\overline
{\lambda}_{i}$ and $(0\leq)$ $\lambda_{j}=\gamma\overline{\lambda}_{j}$, and
so $\gamma=0$ because $\overline{\lambda}_{l}>0$ for all $l\in I$. Hence
$\lambda_{l}=\lambda_{l}^{\prime}$ for all $l\in I$.

Assertion (a) is obvious, while (b) follows immediately from the uniqueness of
the representation of $x\in X$ as $\sum\nolimits_{l\in I}\lambda_{l}u_{l}$
with $\lambda_{l}\geq0$ for $l\in I$ such that $\lambda_{i}=0$ for some $i\in
I$.

(c) Taking $J:=I\setminus\{j\}$ and $J:=I$ in (\ref{r-pg4}) one gets the first
and the second formula, respectively.

(d) Take $x\in X$. Then there exist $(\lambda_{l})_{l\in I}\subseteq
\mathbb{R}$ such that $x=\sum\nolimits_{l\in I}\lambda_{l}u_{l}$. Take
$\alpha:=\min\{\lambda_{l}/\overline{\lambda}_{l}\mid l\in I\}$. Then $\mu
_{l}:=\lambda_{l}-\alpha\overline{\lambda}_{l}\geq0$ for $l\in I$ and $\mu
_{j}=0$ for some $j\in I$, and so $x=\sum\nolimits_{l\in I}\mu_{l}u_{l}\in
Q_{j}$. Therefore, $\bigcup\nolimits_{i\in I}Q_{i}=X=\mathbb{R}^{n}$.
\hfill$\square$

\begin{proposition}
\label{l-GPt4.3}Let $C_{1},C_{2},...,C_{m}\subseteq\mathbb{R}^{n}$ be nonempty
convex sets such that $o\notin C_{i}={}]0,1]C_{i}$ for all $i\in\overline
{1,m}$. Moreover, assume that $\bigcup\nolimits_{j\in J}(\mathbb{P}C_{j}%
)\neq\mathbb{R}^{n}\setminus\{o\}$ for any $J\subseteq\overline{1,m}$ with
$\operatorname*{card}J=n+1$. If $\bigcap\nolimits_{j\in J}C_{j}\neq\emptyset$
for any $J\subseteq\overline{1,m}$ with $\operatorname*{card}J=n$, then
${\bigcap\nolimits_{j\in\overline{1,m}}} C_{j}\neq\emptyset$ and
$\bigcup\nolimits_{j\in\overline{1,m}}(\mathbb{P}C_{j})\neq\mathbb{R}%
^{n}\setminus\{o\}$.
\end{proposition}

Proof. Our aim is to prove that $\bigcap\nolimits_{j\in J}C_{j}\neq\emptyset$
for any $J\subseteq\overline{1,m}$ with $\operatorname*{card}J=n+1$.

Consider $n+1$ sets $A_{1},A_{2},...,A_{n+1}$ among those $m$ given sets.
Using the trick from the proof of Helly's theorem in \cite[p.\ 4]{Sch:93},
take $u_{i}\in A_{1}\cap...\cap\breve{A}_{i}\cap...\cap A_{n+1}$, where
$\breve{A}_{i}$ indicates that $A_{i}$ is the (only) missing set; hence
$u_{i}\in A_{j}$ for all $i,j\in I:=\overline{1,n+1}$ with $i\neq j$, and so
\begin{equation}
\operatorname*{conv}\big\{u_{i}\mid i\in I\setminus\{j\}\big\}\subseteq
A_{j}\quad\forall j\in I. \label{r-pg6}%
\end{equation}

We claim that $o\not \in \operatorname*{conv}\{u_{i}\mid i\in I\}$. Assume
that our claim is true. By Lemma \ref{l-GPt4.1}, after re\-numbering if
necessary, there exist $k\in\overline{1,n}$ and $x\in\operatorname*{conv}%
\{u_{1},...,u_{k}\}\cap\big(]0,1]\operatorname*{conv}\{u_{k+1},...,u_{n+1}%
\}\big)$. As seen above, $u_{i}\in A_{j}$ for all $i,j\in I$ with $i\neq j$,
and so $u_{1},...,u_{k}\in\bigcap\nolimits_{i=k+1}^{n+1}A_{i}$ and
$u_{k+1},...,u_{n+1}\in\bigcap\nolimits_{i=1}^{k}A_{i}$. It follows that
$x\in\bigcap\nolimits_{i=k+1}^{n+1}A_{i}$ and $\alpha^{-1}x\in\bigcap
\nolimits_{i=1}^{k}A_{i}$ for some $\alpha\in{}]0,1]$, that is $x\in
\bigcap\nolimits_{i=1}^{k}\alpha A_{i}\subseteq\bigcap\nolimits_{i=1}^{k}%
A_{i}$. Hence $x\in\bigcap\nolimits_{i\in I}A_{i}$, and so $\bigcap
\nolimits_{i\in I}A_{i}\neq\emptyset$.

Assume now that our claim is not true, that is $o\in\operatorname*{conv}%
\{u_{i}\mid i\in I\}$; then there exists $(\lambda_{i})_{i\in I}%
\subseteq\mathbb{R}_{+}$ such that
\begin{equation}%
%TCIMACRO{\tsum \nolimits_{i\in I}}%
%BeginExpansion
{\textstyle\sum\nolimits_{i\in I}}
%EndExpansion
\lambda_{i}=1,\quad%
%TCIMACRO{\tsum \nolimits_{i\in I}}%
%BeginExpansion
{\textstyle\sum\nolimits_{i\in I}}
%EndExpansion
\lambda_{i}u_{i}=o. \label{r-gp1}%
\end{equation}

(a) Suppose that $\lambda_{j}=0$ for some $j\in I$. Then we get the
contradiction $o=\sum_{i\in I\setminus\{j\}}\lambda_{i}u_{i}\in A_{j}$ by
(\ref{r-pg6}). Hence $\lambda_{i}>0$ for all $i\in I$.

(b) Suppose now that the family $(u_{i})_{i\in I}$ is not a.i.; then there
exists $(\mu_{i})_{i\in I}\subseteq\mathbb{R}$ such that
\begin{equation}%
%TCIMACRO{\tsum \nolimits_{i\in I}}%
%BeginExpansion
{\textstyle\sum\nolimits_{i\in I}}
%EndExpansion
\mu_{i}=0,\quad%
%TCIMACRO{\tsum \nolimits_{i\in I}}%
%BeginExpansion
{\textstyle\sum\nolimits_{i\in I}}
%EndExpansion
\left\vert \mu_{i}\right\vert \neq0,\quad%
%TCIMACRO{\tsum \nolimits_{i\in I}}%
%BeginExpansion
{\textstyle\sum\nolimits_{i\in I}}
%EndExpansion
\mu_{i}u_{i}=o. \label{r-pg3}%
\end{equation}
From (\ref{r-pg3}) we have that $I_{+}:=\{i\in I\mid\mu_{i}>0\}\neq\emptyset$.
Then $\alpha:=\min\{\lambda_{i}/\mu_{i}\mid i\in I_{+}\}>0$, and so $\nu
_{i}:=\lambda_{i}-\alpha\mu_{i}\geq0$ for $i\in I$ and there exists $j\in I$
such that $\nu_{j}=0$. Using (\ref{r-gp1}) and (\ref{r-pg3}) we obtain that
$\sum\nolimits_{i\in I}\nu_{i}=1$ and $\sum\nolimits_{i\in I}\nu_{i}u_{i}=o$,
that is (\ref{r-gp1}) holds with $\lambda_{i}$ replaced by $\nu_{i}$. Using
(a), this is a contradiction because $\nu_{j}=0$. Hence the family
$(u_{i})_{i\in I}$ is affinely independent.

(c) Because $(u_{i})_{i\in I}$ is a.i.\ and $o\in\operatorname*{conv}%
\{u_{i}\mid i\in I\}$ we have that $\mathbb{R}^{n}=\bigcup\nolimits_{i\in
I}Q_{i}$ by Lemma \ref{l-ss}, where $Q_{i}$ is defined in (\ref{r-pg7}).
Clearly, $Q_{j}=\mathbb{R}_{+}\cdot\operatorname*{conv}\big\{u_{i}\mid i\in
I\setminus\{j\}\big\}$ for $j\in I$. Taking into account (\ref{r-pg6}), we
obtain that $\mathbb{R}^{n}=\bigcup\nolimits_{j\in I}Q_{j}\subseteq
\bigcup\nolimits_{j\in I}\mathbb{R}_{+}A_{j}$, and so $\mathbb{R}^{n}%
\setminus\{o\}=\bigcup\nolimits_{j\in I}\mathbb{P}A_{j}$, contradicting our hypothesis.

Therefore, our claim is true, and so $\bigcap\nolimits_{j\in J}C_{j}%
\neq\emptyset$ for any $J\subseteq\overline{1,m}$ with $\operatorname*{card}%
J=n+1$. Because all the sets $C_{j}$ are nonempty, by (the classic) Helly's
theorem (see e.g.\ \cite[Th.\ 1.1.6]{Sch:93}) one obtains that there exists
$\overline{x}\in\bigcap\nolimits_{j\in\overline{1,m}}C_{j}$. Clearly,
$\overline{x}\neq o$, and so $-\overline{x}\notin\mathbb{P}C_{j}$ for all
$j\in\overline{1,m}$, whence $\bigcup\nolimits_{j\in\overline{1,m}}%
(\mathbb{P}C_{j})\neq\mathbb{R}^{n}\setminus\{o\}$. Otherwise, for some
$j^{\prime}\in\overline{1,m}$ and $\lambda\in\mathbb{P}$ one has
$-\lambda\overline{x}\in C_{j^{\prime}}$, and so we get the contradiction
$o=\frac{\lambda}{1+\lambda}\overline{x}+\frac{1}{1+\lambda}(-\lambda
\overline{x})\in C_{j^{\prime}}$. The proof is complete. \hfill\hfill$\square$

\medskip The following result (less $\bigcup\nolimits_{i\in\overline{1,m}%
}K_{i}\neq\mathbb{R}^{n}$) is stated explicitly in \cite[Cor.\ 3]{SarSas:78},
where one refers to \cite[Satz II]{Mol:57} and \cite[Cor.\ p.\ 456]%
{HanKle:69}; in the proof of \cite[Cor.\ 3]{SarSas:78} one uses \cite[Th.\ 2]%
{SarSas:78}.

\begin{corollary}
\label{cor-tHelly-cones}Let $K_{1},K_{2},...,K_{m}\subseteq\mathbb{R}^{n}$,
$m\geq n$, be a family of nontrivial pointed convex cones. If $\bigcap
\nolimits_{j\in J}K_{j}\neq\{o\}$ for any $J\subseteq\overline{1,m}$ with
$\operatorname*{card}J=n$ and $\bigcup\nolimits_{j\in J}K_{j}\neq
\mathbb{R}^{n}$ for any $J\subseteq\overline{1,m}$ with $\operatorname*{card}%
J=n+1$, then $\bigcap\nolimits_{i\in\overline{1,m}}K_{i}\neq\{o\}$ and
$\bigcup\nolimits_{i\in\overline{1,m}}K_{i}\neq\mathbb{R}^{n}$.
\end{corollary}

Proof. Apply Proposition \ref{l-GPt4.3} for $C_{i}:=K_{i}\setminus\{0\}$ for
$i\in\overline{1,m}$. \hfill$\square$

\medskip Lemma \ref{l-ss} shows that the condition \textquotedblleft%
$\bigcup\nolimits_{j\in J}K_{j}\neq\mathbb{R}^{n}$ for any $J\subseteq
\overline{1,m}$ with $\operatorname*{card}J=n+1$" is essential for the
validness of the conclusion of Corollary \ref{cor-tHelly-cones}; similar
remarks are valid for Proposition \ref{l-GPt4.3} and the next result.

\begin{theorem}
\emph{(\cite[Th.\ 4.3]{GuoPen:21}, Helly-type Theorem)} \label{cor-GPt4.3}Let
$C_{1},C_{2},...,C_{m}\subseteq\mathbb{S}$, $m\geq n$, be a family of nonempty
s-convex sets. If $\bigcap\nolimits_{j\in J}C_{j}\neq\emptyset$ for any
$J\subseteq\overline{1,m}$ with $\operatorname*{card}J=n$ and $\bigcup
\nolimits_{j\in J}C_{j}\neq\mathbb{S}$ for any $J\subseteq\overline{1,m}$ with
$\operatorname*{card}J=n+1$, then $\bigcap\nolimits_{i\in I}C_{i}\neq
\emptyset$.
\end{theorem}

Proof. Set $P_{i}:=\mathbb{P}C_{i}$ for $i\in\overline{1,m}$. Because
$\mathbb{P}={}]0,1]\cdot\mathbb{P}$, using Proposition \ref{p1} we have that
$o\notin P_{i}=\operatorname*{conv}P_{i}={}]0,1]P_{i}\neq\emptyset$ for every
$i\in\overline{1,m}$. By hypothesis for any $J\subseteq\overline{1,m}$ with
$\operatorname*{card}J=n$ one has $\emptyset\neq\bigcap\nolimits_{j\in J}%
C_{j}\subseteq\bigcap\nolimits_{j\in J}P_{j}$, and so $\bigcap\nolimits_{j}%
P_{j}\neq\emptyset$. Assume that $\bigcup\nolimits_{j\in J}P_{j}%
=\mathbb{R}^{n}\setminus\{o\}$ for some $J\subseteq\overline{1,m}$ with
$\operatorname*{card}J=n+1;$ then $\mathbb{S}=\rho\big(
\bigcup\nolimits_{j\in J}P_{j}\big)  =\bigcup\nolimits_{j\in J}\rho
(P_{j})=\bigcup\nolimits_{j\in J}C_{j}$, a contradiction. Hence the hypothesis
of Proposition \ref{l-GPt4.3} holds; applying this result we have that there
exists $x\in\bigcap\nolimits_{i=1}^{m}P_{i}$, and so $\rho(x)\in
\bigcap\nolimits_{i=1}^{m}\rho(P_{i})=\bigcap\nolimits_{i=1}^{m}C_{i}$.
\hfill\hfill$\square$

\medskip Having in view the equivalence (i)~$\Leftrightarrow$ (ix) from
Proposition \ref{p1}, in the preceding proof we could take $]0,1]C_{i}$
instead of $\mathbb{P}C_{i}$ for $i\in\overline{1,m}$ when $\Phi$ is,
moreover, convex (hence sublinear).

\begin{theorem}
\emph{(\cite[Th.\ 4.8]{GuoPen:21}, Carath\'{e}odory-type Theorem)}
\label{cor-GPt4.8}Let $S\subseteq\mathbb{S}$ be a nonempty hull-addible set.
Then for each $x\in\operatorname*{sco}S$, there exist $(x_{i})_{i=1}%
^{n}\subseteq S$ and $(\lambda_{i})_{i=1}^{n}\subseteq\mathbb{P}$ with
$\sum_{i=1}^{n}\lambda_{i}=1$ such that $x=$\emph{(s)}$\sum_{i=1}^{n}%
\lambda_{i}x_{i}$.
\end{theorem}

Proof. Take $x\in\operatorname*{sco}S$. Then $x\in\mathbb{R}_{+}%
\operatorname*{sco}S=\mathbb{R}_{+}\rho(\operatorname*{conv}S)=\mathbb{R}%
_{+}\operatorname*{conv}S=\operatorname*{conv}\left(  \mathbb{R}_{+}S\right)
$. Because $x\neq o$, using Lemma \ref{l-pg1}~(ii), there exist $(x_{i}%
)_{i=1}^{n}\subseteq S$ and $(\lambda_{i})_{i=1}^{n}\subseteq\mathbb{R}_{+}$
such that $x=\sum_{i=1}^{n}\lambda_{i}x_{i}$; then $\lambda:=\sum_{i=1}%
^{n}\lambda_{i}>0$. Setting $\lambda_{i}^{\prime}:=\lambda^{-1}\lambda_{i}$
$(\geq0)$, we have that $\sum_{i=1}^{n}\lambda_{i}^{\prime}=1$ and so
$\lambda^{-1}x=:x^{\prime}=\sum_{i=1}^{n}\lambda_{i}^{\prime}x_{i}$. Hence
$x=\rho(x^{\prime})=$(s)$\sum_{i=1}^{n}\lambda_{i}^{\prime}x_{i}$.
\hfill$\square$

\medskip Consider $C\subseteq\mathbb{S}$ an s-convex set. One says that $x\in
C$ is an \emph{extreme point} of $C\subseteq\mathbb{S}$ if $[x_{1},x_{2}\in
C$~$\wedge$ $\lambda\in{}]0,1[$~$\wedge$ $x=\lambda x_{1}+_{s}(1-\lambda
)x_{2}]$~$\Rightarrow$ $x_{1}=x_{2}$; the set of s-extreme points of $C$ is
denoted by $\operatorname*{sext}C$.

\begin{proposition}
\label{l-p5.2}\emph{(\cite[Prop.\ 5.3 and Cor.\ 5.5]{GuoPen:21}) }Let
$C\subseteq\mathbb{S}$ be s-convex and $x\in C$. Then

\emph{(i)}~$x\in\operatorname*{sext}C$ if and only if $\mathbb{R}_{+}x$ is an
extreme ray of $\mathbb{R}_{+}C$.

\emph{(ii)}~Assume that $x:={}$\emph{(s)}$\sum_{i=1}^{k}\lambda_{i}x_{i}%
\in\operatorname*{sext}C$ with $(x_{i})_{i=1}^{k}\subseteq C$, $(\lambda
_{i})_{i=1}^{k}\subseteq\mathbb{P}$ and $\sum_{i=1}^{k}\lambda_{i}=1;$ if
$x\in\operatorname*{sext}C$ then $x_{i}=x$ for $i\in\overline{1,k}$.
\end{proposition}

Proof. Notice that having the pointed convex cone $P\subseteq\mathbb{R}^{n}$,
$\mathbb{R}_{+}x$ with $x\in P\setminus\{0\}$ is an extreme ray (that is, a
face) of $P$ if and only if for all $y,z\in P\setminus\{0\}$ one has
$x\in\mathbb{P}y+\mathbb{P}z~\Rightarrow$ $\mathbb{P}x=\mathbb{P}%
y=\mathbb{P}z$. By induction one obtains that $\mathbb{R}_{+}x$ with $x\in
P\setminus\{0\}$ is an extreme ray of $P$ if and only if for all
$x_{1},...,x_{k}\in P\setminus\{0\}$ $(k\geq2)$ one has $x\in\mathbb{P}%
x_{1}+...+\mathbb{P}x_{k}~\Rightarrow$ $\mathbb{P}x_{i}=\mathbb{P}x$ for
$i\in\overline{1,k}$.

Using Proposition \ref{p1} we have that $K:=\mathbb{R}_{+}C$ is a pointed
convex cone.

(i) Assume that $x\in\operatorname*{sext}C$ $(\subseteq K\setminus\{0\})$ and
take $y,z\in K\setminus\{0\}$ such that $x\in\mathbb{P}y+\mathbb{P}z$; setting
$y^{\prime}:=\rho(y)$, $z^{\prime}:=\rho(z)$, we have that $y^{\prime
},z^{\prime}\in C$ and $x\in\mathbb{P}y^{\prime}+\mathbb{P}z^{\prime}$. Hence
$x=\alpha y^{\prime}+\beta z^{\prime}$ with $\alpha,\beta\in\mathbb{P}$,
whence $\frac{1}{\alpha+\beta}x=\frac{\alpha}{\alpha+\beta}y^{\prime}%
+\frac{\beta}{\alpha+\beta}z^{\prime}$. It follows that $x=\rho\big(\frac
{1}{\alpha+\beta}x\big)=\frac{\alpha}{\alpha+\beta}y^{\prime}+_{s}\frac{\beta
}{\alpha+\beta}z^{\prime}$. Because $\frac{\alpha}{\alpha+\beta}\in{}]0,1[$
and $x\in\operatorname*{sext}C$, we obtain that $y^{\prime}=z^{\prime}$,
whence $\mathbb{P}y^{\prime}=\mathbb{P}z^{\prime}$. Hence $\mathbb{P}%
y=\mathbb{P}z$, which proves that $\mathbb{R}_{+}x$ is an extreme ray of $K$.

Conversely, assume that $\mathbb{R}_{+}x$ is an extreme ray of $K$ and take
$y,z\in C$, $\lambda\in{}]0,1[$ such that $x=\lambda y+_{s}(1-\lambda)z$. It
follows that there exists $\gamma\in\mathbb{P}$ such that $\gamma x=\lambda
y+(1-\lambda)z$, whence $x\in\mathbb{P}y+\mathbb{P}z$; because $\mathbb{R}%
_{+}x$ is an extreme ray of $\mathbb{R}_{+}C$ we have that $\mathbb{R}%
_{+}y=\mathbb{R}_{+}z$, and so $\{y\}=\mathbb{S}\cap\mathbb{R}_{+}%
y=\mathbb{S}\cap\mathbb{R}_{+}z=\{z\}$. Therefore, $x\in\operatorname*{sext}C$.

(ii) Clearly, $\gamma x=\sum_{i=1}^{k}\lambda_{i}x_{i}$ for some $\gamma>0$,
and so $x\in\mathbb{P}x_{1}+...+\mathbb{P}x_{k}$. Because $\mathbb{R}_{+}x$ is
an extreme ray of $\mathbb{R}_{+}C$, we have that $\mathbb{P}x=\mathbb{P}%
x_{i}$ for $i\in\overline{1,k}$, whence $x=x_{i}$ for $i\in\overline{1,k}$.
\hfill$\square$

\medskip Note that assertion (i) of Proposition \ref{l-p5.2} is
\cite[Prop.\ 5.3]{GuoPen:21}, while assertion (ii) is established in
\cite[Cor.\ 5.5]{GuoPen:21} under the additional assumption that $C$ is closed.

\begin{theorem}
\label{p-GPt5.5}\emph{(\cite[Th.\ 5.6 and Cor.\ 5.4]{GuoPen:21})} Let
$C\subseteq\mathbb{S}$ be closed and s-convex. Then $C=\operatorname*{sco}%
(\operatorname*{sext}C);$ in particular, $\operatorname*{sext}C\neq\emptyset$.
\end{theorem}

Proof. Because $\mathbb{R}_{+}C$ is a pointed closed convex cone, using
\cite[Cor.\ 18.5.2]{Roc:70} and Proposition \ref{l-p5.2}~(i) we have that
$\mathbb{R}_{+}C=\operatorname*{cone}(\operatorname*{sext}C)$; in particular
$\operatorname*{sext}C\neq\emptyset$.

Of course, $\operatorname*{sco}(\operatorname*{sext}C)\subseteq C$. Take $x\in
C$ $[\subseteq(\mathbb{R}_{+}C)\setminus\{o\}]$; hence there exist
$(x_{i})_{i=1}^{k}\subseteq\operatorname*{sext}C$ and $(\lambda_{i})_{i=1}%
^{k}\subseteq\mathbb{P}$ such that $x=\sum_{i=1}^{k}\lambda_{i}x_{i}$. Then
\[
x=\rho(x)=\rho\big(%
%TCIMACRO{\tsum \nolimits_{i=1}^{k}}%
%BeginExpansion
{\textstyle\sum\nolimits_{i=1}^{k}}
%EndExpansion
\lambda_{i}x_{i}\big)=\rho\big(%
%TCIMACRO{\tsum \nolimits_{i=1}^{k}}%
%BeginExpansion
{\textstyle\sum\nolimits_{i=1}^{k}}
%EndExpansion
\lambda_{i}^{\prime}x_{i}\big)=\rho\big(%
%TCIMACRO{\tsum \nolimits_{i=1}^{k}}%
%BeginExpansion
{\textstyle\sum\nolimits_{i=1}^{k}}
%EndExpansion
\lambda_{i}^{\prime}x_{i}\big)=(\text{s})%
%TCIMACRO{\tsum \nolimits_{i=1}^{k}}%
%BeginExpansion
{\textstyle\sum\nolimits_{i=1}^{k}}
%EndExpansion
\lambda_{i}^{\prime}x_{i},
\]
where $\lambda_{i}^{\prime}:=\lambda_{i}/\sum_{i=1}^{k}\lambda_{i}^{\prime}$
for $i\in\overline{1,k}$. Therefore, $C=\operatorname*{sco}%
(\operatorname*{sext}C)$. \hfill$\square$

\end{document}